\documentclass{article}

 \usepackage{graphicx} 
 \usepackage{float}
 \usepackage{amsmath} 
\usepackage[utf8]{inputenc} % allow utf-8 input
\usepackage[T1]{fontenc}    % use 8-bit T1 fonts
\usepackage{hyperref}       % hyperlinks
\usepackage{url}            % simple URL typesetting
\usepackage{booktabs}       % professional-quality tables
\usepackage{amsfonts}       % blackboard math symbols
\usepackage{nicefrac}       % compact symbols for 1/2, etc.
\usepackage{microtype}      % microtypography
\usepackage{xcolor}         % colors
\usepackage[nonatbib]{neurips_2024}
\usepackage{algorithm, algorithmic}
\title{Stochastic Optimization Using Ricci Flow}

\author{%
  Varsha Gupta\thanks{Use footnote for providing further information
    about author (webpage, alternative address)---\emph{not} for acknowledging
    funding agencies.} \\
  Department of Agricultural and Biological Engineering\\
  Purdue University\\
  West Lafayette, IN 47906 \\
  \texttt{vvarsha@purdue.edu} \\
}

\begin{document}

\maketitle

\begin{abstract}
This paper proposes a theoretical framework for modeling and optimizing the bounded functions based on the Fourier series approximation and Ricci flow. Specifically, the initial manifold, $\mathcal{M}_0$ is approximated using Fourier series approximation in conjunction with the center and boundary sampling procedure introduced in the paper. The manifold is iteratively evolved using an algorithm that involves sampling along geodesic hyper-sphere defined by the Riemannian metric tensor. Thus obtained surrogate manifold is optimized by applying inverse Ricci flow i.e. instead of regularizing the manifold, flow allows for the high curvature regions to blow into finite time singularities. This allows for the singularities to occur at potential global optima assuming the deviation of the manifold at any point is smaller than the optimum. In addition, the error bound is established on the accuracy of the surrogate manifold. Finally, the proposed method is tested on stochastic sampling from five benchmark functions to illustrate the utility of this method. 
\end{abstract}

\section{Introduction}
Due to intrinsic uncertainty, stochastic processes are far more challenging to model and optimize in comparison to the corresponding deterministic processes[1]. From its inception [2,3], a significant body of work dedicated to modeling and optimization of these processes has since been developed [4,5].\\
The modeling approach includes approximating the function to be log-concave [6,7]. Therefore, the results obtained using these methods are limited in their ability to predict the local optimum of highly non-convex functions.\\
Moreover, a large amount of the work skews towards using unidirectional randomized search methods such as random walk [8,9] and hit-and-run [10], and simulated annealing [11]. In addition, zeroth order optimization [12] is a favored method to reduce computational complexity [13]. \\
Monte Carlo-based sampling methods for optimizing stochastic optimization find applications in a variety of problems such as machine learning, supply chain networking, engineering design, and scheduling [14]. However, these methods are effective only when a small sample size is needed. Owing to the curse of dimensionality, stochastic optimization suffers from large computation complexity as the sample size increases drastically [15]. \\
To address the problem of optimizing highly non-convex functions with reasonable accuracy, we propose a theoretical framework outlined in the paper. The motivation of this work is to exploit the topological structure embedded in the data and use appropriate methods to deform this structure strategically toward obtaining optimum solutions.\\
The paper is broadly divided into two sections. In section 2, wave superposition is used to approximate the function manifold. A detailed sampling procedure is outlined to approximate the initial manifold and its evolution. In section 3, the manifold is optimized using a strategic combination of different types of Ricci flows. The flow allows for singularities to happen at optimum locations. These locations are filtered based on the sign of the objective function (maxima or minima). The filtered points are checked and compared for the best feasible solution.\\
The method outlined in the paper is a theoretical framework that can serve as a basis for solving highly non-convex problems. However, the method can be tailored for specific problems depending on the problem specifications. Section 4 is dedicated towards establishing error bounds on surrogate manifold accuracy compared to the original function. Then we apply our method to stochastic samples collected from five different benchmark function in section 5. Here, we discuss the results in detail especially for highly non-convex function such as Rastrigin. Finally, we close-off with conclusion of our paper briefly outlined in section 6.

\section{Wave superposition for topology approximation}
The motivation for modeling the stochastic process is to extract the topology embedded in the sampled dataset. The superposition of sine waves or Fourier series approximation can approximate a manifold using a set of samples. The manifold approximation allows for the sampling function to become smooth and continuous. In other words, the information is extrapolated in the regions beyond the sampling point in a smooth continuous space in addition to becoming infinitely differentiable. This allows us to use differential geometry based sophisticated optimization techniques that search of the optimum using the curvature information.

\textbf{Assumptions:}
The following assumptions are made in order to approximate and optimize the function using the proposed method: 
\begin{itemize}
    \item The function can be approximated reasonably accurately such that deviation at no point is larger than the optimum value of the function.
    \item The manifold is smooth.
    \item The function is bounded and continuous.  
    \end{itemize}
\subsection{Manifold approximation}
Let $f:\mathbb{R}^n \to \mathbb{R}$ be a bounded sampling function such that 
\[ 
f: \prod_{d \in \mathcal{D}}[d_L, d_U] \to \mathbb{R}, \,\,\,\,\,\,\,\,
\mathcal{D} = \{1,2,3,...,n\}
\]
The function can be approximated at any point $\mathbf{x}$ in the sampled domain using the Fourier series transformation as follows:
\begin{equation}\label{FT}
    f(\mathbf{x}) \approx \sum_{\mathbf{k}\in K} a_{\mathbf{k}} e^{j \mathbf{k} \omega \mathbf{x}}
\end{equation}
where $\omega$ represents fundamental frequency in each dimension and K is a set of indices of all the frequency components.
We will denote the R.H.S. in the above equation as $F(x)$. 

The coefficients $a_k$ can be approximated using $N$ samples as follows:
\begin{equation}\label{coeff}
    a_{\mathbf{k}} = \frac{1}{N} \sum_{i=1}^N f(\mathbf{x}_i) e^{-j \mathbf{k} \omega \mathbf{x}_i} 
\end{equation}
If we define the error function as following:
\begin{equation}\label{error}
 \epsilon (N) = |f(x)-F(x)|
\end{equation}
Then the error function approaches zero as the number of samples reaches to infinity, i.e.
\[
\mbox{As } N\to \infty,\,\,\  \epsilon \to 0 
\]
Please note that due to the superposition of sine waves, the manifold is smooth in nature. The number of waves superimposed in the same direction favoring the optima provides local information about the distribution of optimum points. However, in order to minimize sampling, we adopt a more detailed procedure for approximating the manifold as outlined below.

\subsection{Sampling procedure and manifold evolution}
We begin by approximating the manifold using a minimal number of sampled data points that allow for the manifold to best approximate the functions at the corner points and center point. 

Therefore, for the first iteration, we define sampling using the following two sets of samples:

\begin{enumerate}
    \item Boundary sampling: These samples are ideally collected corresponding to the corner points in all the dimensions. These points can alternately be taken as the extreme points in the feasible space. However, if the feasible space is unknown, these points can be taken arbitrarily such that it allows for a large set of data points to lie in this region. These sets are chosen based on the specifics of the problem under consideration. 
    
    \item Midpoint sampling: The points are ideally evaluated at the midpoint such that each coordinate corresponds to the midpoint of the respective dimensions bounds. We denote this midpoint as $\mathbf{p}$. 
\end{enumerate}

Using these samples and by adopting the outlined procedure for approximating the manifold in section 2.1, we approximate the manifold. Let's call this surrogate manifold, $\mathcal{M}_0^{n+1}$

Next, we draw circles on this manifold as it evolves at constant intervals such that the radii of each of these circles,$C_{\mathbf{p},r}$ on the z-evolved manifold, $\mathcal{M}_z^{n+1}$ are given by the Riemannian metric tensor from the midpoint to the surface of the sphere. 

\[
C_{\mathbf{p},r} = \{\mathbf{q} \in \mathcal{M}_z^{n+1}|\,\,\,\ d(p,q) = r\}
\]

where g is geodesic distance from $\mathbf{p}$ to $\mathbf{q}$. For dimensions greater than 2, we have geodesic hyper-sphere instead of circles.

For each circle, we randomly pick a point and evaluate the function, F(x) value at this point. 
\[
\mbox{if  } |f(x) - F(x) \leq \alpha, \mbox{we reject the sample}\]
\[
\mbox{ else,  we accept the sample.}
\]

At each iteration, we collect the samples and refine the surrogate manifold by using the Fourier series approximation method in section 2.1.

We continue the procedure until the sampling circle $C_{\mathbf{p},r}$ extends beyond the manifold's domain boundary.
The algorithm is written as follows.

\begin{algorithm}
\caption{Manifold Approximation and Sampling}
\begin{algorithmic}[1]
\STATE Initialize the number of dimensions, $n$.
\STATE Define the bounds for each dimension.

\STATE \textbf{Boundary Sampling:}
\STATE Collect samples at the corner points in all dimensions.
\STATE If the feasible space is unknown, choose extreme points arbitrarily.

\STATE \textbf{Midpoint Sampling:}
\STATE Sample the midpoint, $\mathbf{p}$.

\STATE \textbf{Approximate Initial Manifold:}
    \STATE Approximate the initial surrogate manifold, denoted as $\mathcal{M}_0^{n+1}$.
    
    \STATE \textbf{Circular sampling and surrrogate Manifold evolution}
    \FOR{each iteration}
        \STATE \textbf{Sampling on Circles:}
         \STATE Draw circles or geodesic hyper-spheres centered at $\mathbf{p}$ with radii determined by the Riemannian metric tensor.
            \FOR{each circle $C_{\mathbf{p},r}$ defined by $d(\mathbf{p}, \mathbf{q}) = r$}
                \STATE Randomly select a point $\mathbf{q}$ on the circle.
                \STATE Evaluate the function $F(\mathbf{q})$ at point $\mathbf{q}$.
            \ENDFOR
    
            \STATE \textbf{Sample Acceptance Criteria:}
            \IF{$|f(\mathbf{q}) - F(\mathbf{q})| \leq \alpha$}
                \STATE Reject the sample.
            \ELSE
                \STATE Accept the sample.
            \ENDIF
   
    \STATE Evolve the surrogate manifold $\mathcal{M}_z^{n+1}$ at constant intervals using accepted samples.
\ENDFOR

\STATE \textbf{Check Termination:}
\STATE Continue the iterations until the circle $C_{\mathbf{p},r}$ crosses the boundary of the manifold.

\STATE \textbf{Output:}
\STATE The final surrogate manifold.

\end{algorithmic}
\end{algorithm}

The final surrogate manifold is then optimized using appropriate flow conditions on the manifold. The flow should be such that it allows the detection of extreme points on the manifold. 

\section{Ricci Flow Informed Optimization}
We construct a Riemannian metric $g(t)$ on the manifold $\mathcal{M}$ so that large values of the objective function $f(\mathbf{x})$ correspond to high curvature regions. A typical choice is 
\[
  g_{ij}(\mathbf{x}) \;=\; e^{\,\beta f(\mathbf{x})} \,\delta_{ij},
\]
where $\beta > 0$ is a scaling constant and $\delta_{ij}$ is the Kronecker delta function. Under the inverse Ricci flow,
\begin{equation}
\label{eq:invRicciFlow}
  \frac{\partial g}{\partial t} \;=\; +2\,\mathrm{Ric}(g),
\end{equation}
regions already exhibiting high curvature tend to blow-up. Classical results in geometric analysis [16, 17] show that these high-curvature peaks can blow-up to finite-time singularities. Hence, if $f$ has a global maximum at $\mathbf{x}^*$, the metric curvature around $\mathbf{x}^*$ tends to diverge, and the flow forms a singularity precisely at these locations. In practice, we track these singularities as candidate global optima of $f$. \\
To find the optimum on the manifold, we apply an iterative process using a combination of Ricci flow and inverse Ricci flow as defined below.

In the first iteration, for the Riemannian metric, g(t), the manifold, $\mathcal{M}_z^{n+1}$ for t $\in (0,T]$ is deformed using inverse Ricci flow which is defined as follows:

\[
\frac{\partial{g}}{\partial{t}} = 2*RiC(g)
\]
This inverse Ricci flow causes the space to deform such that the point of largest curvature increases until the surface tends to blow-up at this point into singularity.\\

In the next iteration, we apply the inverse Ricci flow on  $\mathcal{M}_z^{n+1}$ with the exception that at the point approaching singularity in the previous iteration:

\[
\frac{\partial{g}}{\partial{t}} = -2*RiC(g)
\]

This process continues until the desired number of candidate points are gathered for optima. Please note that the application of Ricci flow here helps regularize the surrogate manifold as alternate to surgery thereby allowing the optimization process to continue. \\

We approximated the manifold using Fourier series approximation. The candidate points, filtered according to whether the objective is to minimize or maximize, are evaluated to select the best solution.

\begin{algorithm}[H]
\caption{Finding the Optimum on the Manifold}
\begin{algorithmic}[2]
\STATE Initialize the manifold $\mathcal{M}_z^{n+1}$ and the Riemannian metric $g(t)$ for $t \in (0, T]$.

\STATE \textbf{First Iteration: Inverse Ricci Flow}
\FOR{each $t \in (0, T]$}
    \STATE Deform the manifold $\mathcal{M}_z^{n+1}$ using the inverse Ricci flow:
    \[
    \frac{\partial g}{\partial t} = 2 \cdot \text{Ric}(g)
    \]

\ENDFOR

\STATE \textbf{Next Iteration: Ricci Flow}
        \IF{curvature tends to blow from the previous iteration}
                \STATE Deform the manifold $\mathcal{M}_z^{n+1}$ using the Ricci flow:
        \[
        \frac{\partial g}{\partial t} = -2 \cdot \text{Ric}(g)
        \]
        \ELSE
                \STATE Deform the manifold $\mathcal{M}_z^{n+1}$ using the inverse Ricci flow:
    \[
    \frac{\partial g}{\partial t} = 2 \cdot \text{Ric}(g)
    \]
            \ENDIF
\STATE Continue until the desired number of candidate points for optima are reached.

\STATE \textbf{Filtering Candidate Points}
\STATE Filter the candidate points based on the objective function (minima or maxima).

\STATE Choose the best solution among the filtered candidate points as the local optimum.

\STATE \textbf{Output:}
\STATE The local optimum on the manifold.

\end{algorithmic}
\end{algorithm}

\section{Error Bound on Surrogate Manifold}

We know that when a Riemannian metric is defined such that optima of $f(x)$ lead to regions of high curvature, applying the inverse Ricci flow tends to amplify these curvature peaks, ultimately resulting in finite-time singularities at the global maximum of the function [18].\\
Therefore, the stochastic optimization process using the proposed method is limited by Algorithm~1's ability to approximate the original function. Consequently, in order to establish bounds on this approach, we shall derive the error bounds for the Fourier-based approximation presented in the paper. \\

We know that the error in approximating the original function arises on account of two factors:

\begin{itemize}
    \item Fourier series approximation using truncated series, $(E_F)$.
    \item Finite stochastic sampling error to estimate Fourier coefficients, $(E_S)$.
\end{itemize}

If $f$ is in a Sobolev space of order $s$ and is sufficiently smooth, then its Fourier coefficients shall decay quickly. Truncating the series to $N$ terms then gives an approximation error bounded by following [19]:

\begin{equation}\label{eq:fourier_error}
E_F = |f - F_N| \leq C_F\,N^{-s},
\end{equation}

where $C_F$ is a function of smoothness of $f$ and on the chosen norm basis.

Now, for deterministic cases, the sampling error can be only due to finite number of samples. However, for stochastic sampling, we need to account for the intrinsic noise.\\
First, we establish an error bound for finite sampling. Because we only sample $f$ at a finite set of points, the sampling error arises. Let's say, $f$ satisfies a Lipschitz condition,
\begin{equation}
| f(x) - f(y) | \;\leq\; L\,|x-y|,
\end{equation}

and if the maximum gap among neighboring samples in the domain is $\Delta$ then we can bound sampling error as following:

\begin{equation}\label{eq:recon_error}
E_S \leq C_S \, \Delta,
\end{equation}

where $ C_S $ is a constant.\\

In a uniform grid of dimension $n$, we can simplify the above expression using $\Delta \sim N^{-1/n}$ as following: 
\begin{equation}\label{eq:recon_error}
E_S \leq C_S \, N^{-1/n},
\end{equation}

In a stochastic setting, we need to account for noise. For large enough samples, we can refine the above bound as following [20]:
\begin{equation}
E_S \leq C_S \,  N^{-1/n} + C_{\sigma} \,\sqrt{\frac{\log(1/\delta)}{N}},
\end{equation}

where $ C_{\sigma} $ is a constant, $N$ is the total number of samples, and $\delta$ is the confidence interval. 

Therefore, the total error can be written in the form of following upper bound:

\begin{equation}\label{eq:combined_error}
E_{\text{total}} \leq C_F\,N^{-s} + C_S\,N^{-1/n} + C_{\sigma} \,\sqrt{\frac{\log(1/\delta)}{N}},
\end{equation}

Based on this error bound, we can say that:

\begin{itemize}
  \item \textbf{Fourier Error ($E_F$)} decreases as we include higher order terms. For smooth functions, a close approximation of the original function is possible with lower-order terms.
  \item \textbf{Sampling Error ($E_S$)} decreases as the sampling density increases i.e. as $\Delta \to 0$ and $N \to \infty$. 
\end{itemize}

Please note that for smooth functions, the first term will be small. Conversely, for large number of terms needed to approximate the oscillatory functions, this term will be large. Similarly, the second term can be small or large depending on the number of dimensions. However, noise $\sigma$ remains a limiting factor in stochastic environments, although it can be reduced by averaging with a large number of samples.

\section{Results}
To test the algorithm outlined in the paper, we create randomized sampling from deterministic benchmark functions to simulate stochastic processes. In the first step, we develop the surrogate manifold within the domain specified for each function using the Algorithm~1. Once the manifold is developed, we apply Algorithm~2 to get optimal point. \\

We select five benchmark functions: Rosenbrock function [21], Himmelblau [22], Booth function [23], Ackley function [24] and Rastrigin function [25]. To allow for stochasticity, we randomly generate samples from each of these functions. Coding is done in python using Jupyter Notebook.

\subsection{Manifold Approximation}
\begin{figure}[htbp!]
    \hspace*{-2cm} 
    \includegraphics[width=1.2\linewidth]{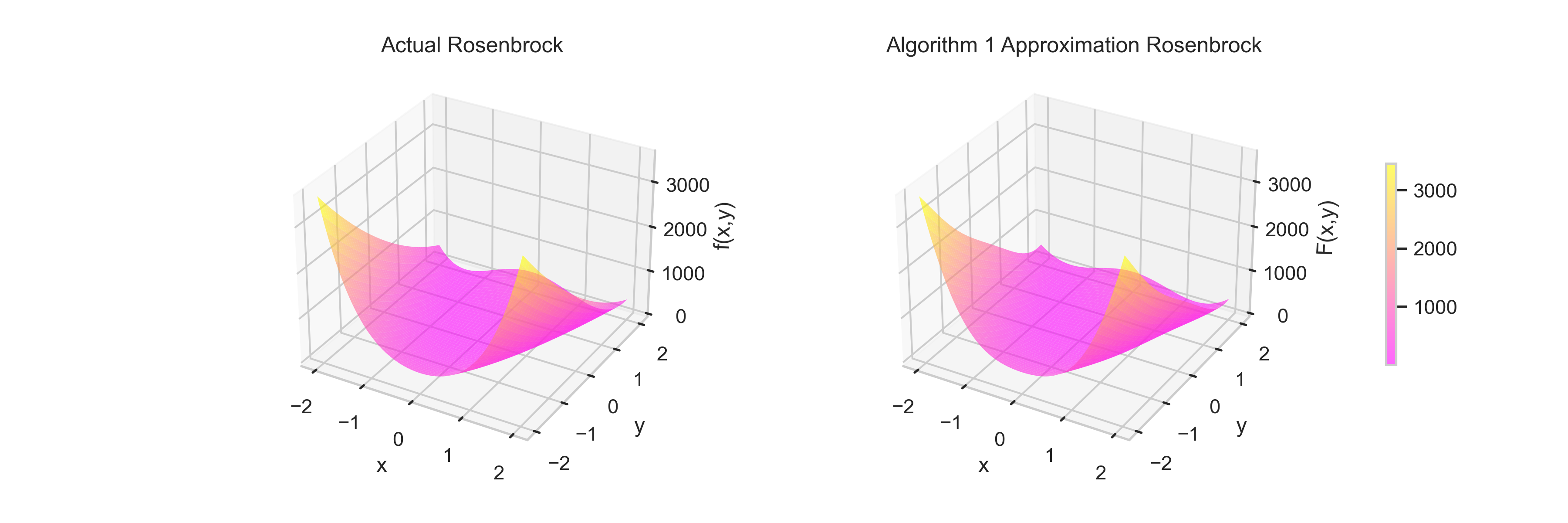}
    \caption{Rosenbrock approximation using stochastic sampling}
    \label{Rosenbrock}
\end{figure}

\begin{figure}[htbp!]
    \hspace*{-2cm} 
    \includegraphics[width=1.2\linewidth]{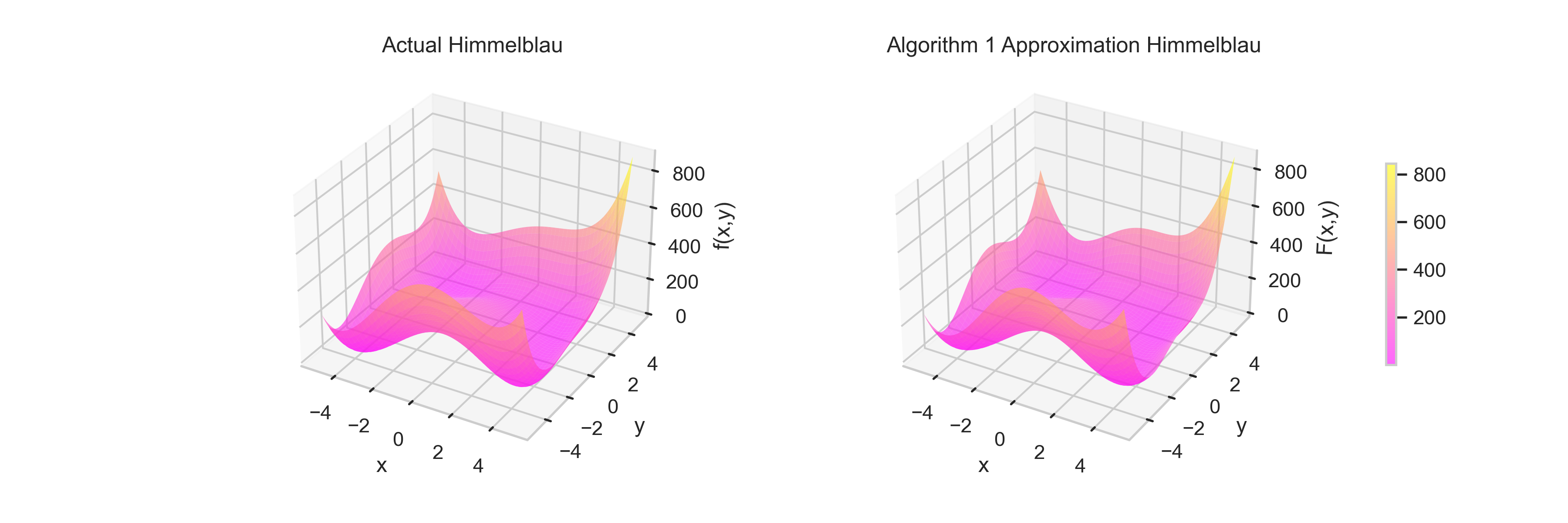}
    \caption{Himmelblau approximation using stochastic sampling}
     \label{Himmelblau}
\end{figure}

\begin{figure}[htbp!]
    \hspace*{-2cm} 
    \includegraphics[width=1.2\linewidth]{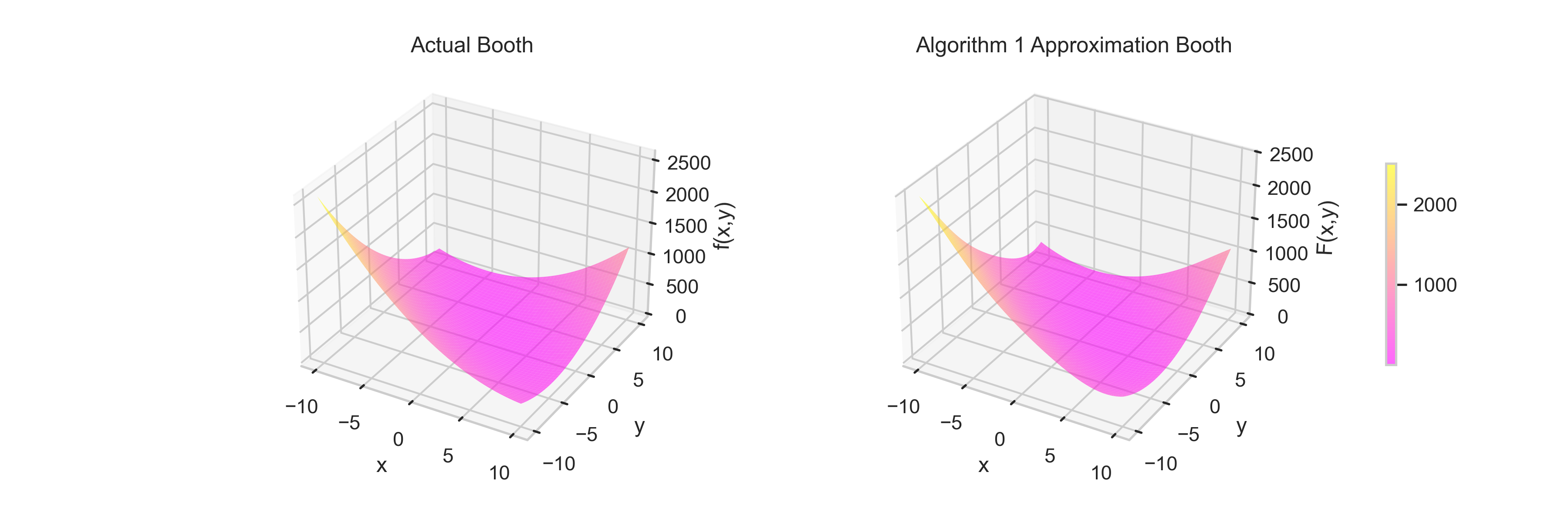}
    \caption{Booth approximation using stochastic sampling}
     \label{Booth}
\end{figure}

\begin{figure}[htbp!]
    \hspace*{-2cm} 
    \includegraphics[width=1.2\linewidth]{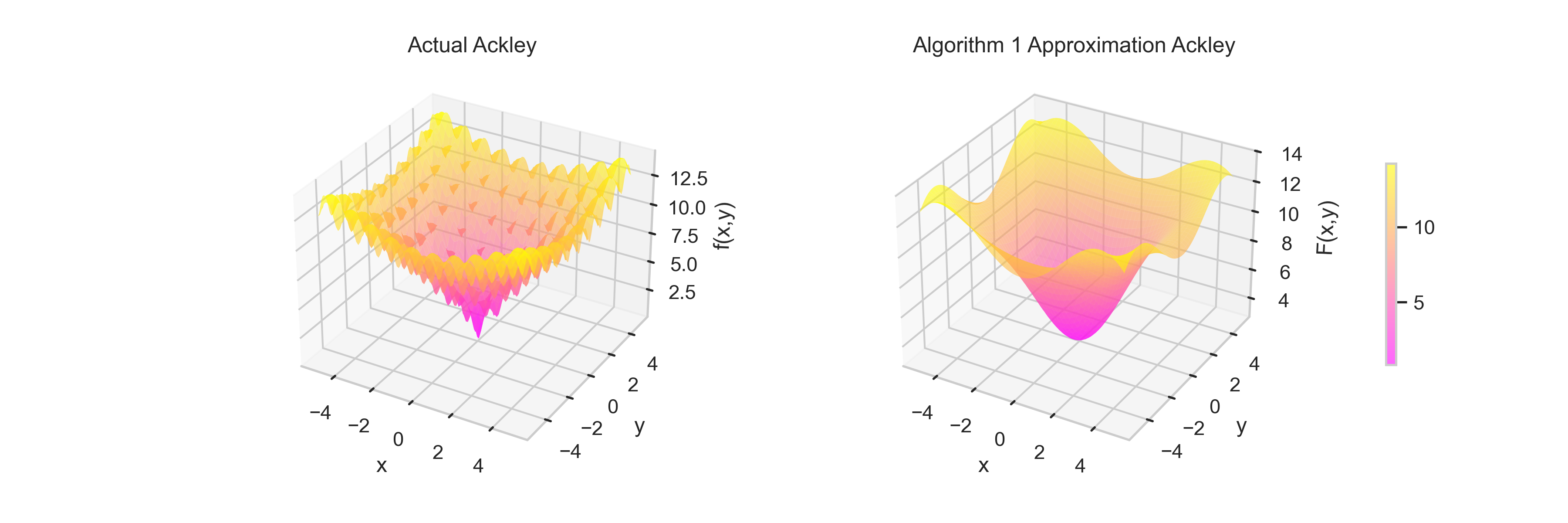}
    \caption{Ackley approximation using stochastic sampling}
     \label{Ackley}
\end{figure}

The Fourier series approximation is based on a Fourier expansion of order $M=3$, a sample acceptance threshold of $\alpha=0.1$, and circular sampling with radii chosen based on each function's domain.  The resulting error metrics for each benchmark function is summarized in Table~\ref{tab:fourier}. In addition, visual comparison for each function's approximation against the respective actual function  are shown in Fig. ~\ref{Rosenbrock},  Fig.~\ref{Himmelblau}, Fig. ~\ref{Booth},  Fig.~\ref{Ackley}, Fig. ~\ref{Rastrigin}.\\

\begin{table}[htbp!]
\centering
\caption{Algorithm~1 Results}
\label{tab:fourier}
\begin{tabular}{l l l l l l l}
\toprule
\textbf{Function} & \textbf{Domain}  & \textbf{MAE} & \textbf{MSE} & \textbf{R$^2$} \\
\midrule
Rosenbrock  & $[-2,2]\times[-2,2]$        & 17.7524  & 1041.8165 & 0.9974 \\
Himmelblau & $[-5,5]\times[-5,5]$        & 7.3665   & 137.8512  & 0.9899 \\
Booth       & $[-10,10]\times[-10,10]$     & 7.2568   & 202.1116  & 0.9990 \\
Ackley      & $[-5,5]\times[-5,5]$        & 0.5261   & 0.4350    & 0.9324 \\
Rastrigin   & $[-5.12,5.12]\times[-5.12,5.12]$  & 8.3745   & 106.2288  & 0.4887 \\
\bottomrule
\end{tabular}
\end{table}

\begin{figure}[htbp!]
   \hspace*{-2cm} 
    \includegraphics[width=1.2\linewidth]{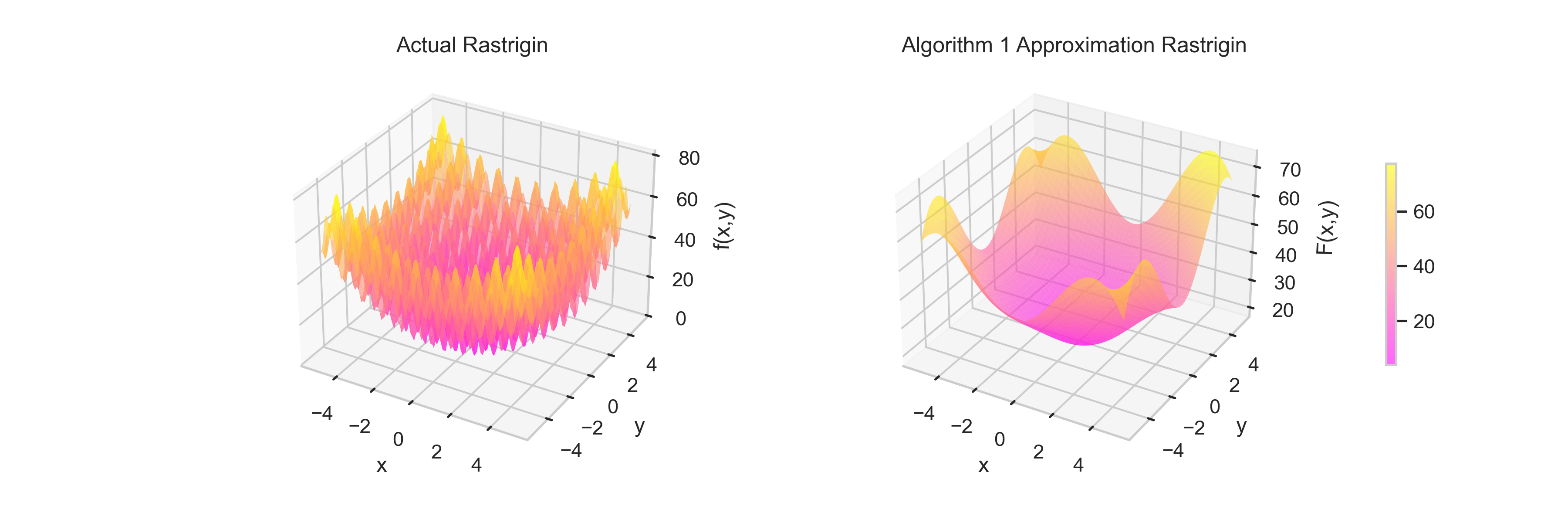}
    \caption{Rastrigin approximation using stochastic sampling}
     \label{Rastrigin}
\end{figure}

Following are the major insights obtained from the approximation of surrogate manifolds:

\begin{itemize}
    \item The surrogate manifolds are infinitely differentiable due to sine wave superposition. This is crucial for Ricci flow-based optimization in Algorithm~2 to exploit curvature for finding optima.

    \item Even with low-order Fourier approximation, relatively smoother functions such as Himmelblau and Booth functions visually resemble the standard function shape. The error metrics also confirm these results.

    \item Rosenbrock has a very strong $R^2$ $\approx$ 0.99 but high MAE and MSE values indicate that although the model over-estimates absolute errors in some regions, it approximates the overall shape well.
    
    \item Oscillatory functions such as Ackley and Rastrigin show deviation from the standard shape especially near local optimums.  As Table~\ref{tab:fourier} shows, the truncated Fourier representation achieves a relatively lower values of $R^2$ $\approx$ 0.93 for Ackley and $R^2$ $\approx$ 0.49 for Rastrigin. However, the iterative sampling, based on geodesic circles in the manifold, still refines the approximation enough to retain a qualitatively correct global structure of peaks and valleys.\\
\end{itemize}

\subsection{Surrogate Manifold Optimization}
Now, that we have a surrogate manifold, we shall apply Ricci flow based Algorithm~2 to find the optimal. Following are the common simulation parameters across each function:
\begin{itemize}
  \item Time step: $dt = 0.001$
  \item Number of iterations: $300$
  \item Curvature threshold: $1 \times 10^{-5}$
  \item Grid resolution: $200 \times 200$
  \item Plot interval: $300$
\end{itemize}

To evaluate Algorithm~2 on purely deterministic functions, we first apply it to locate the optima of each benchmark function in a noise-free setting. Table~\ref{tab:opt} compares the candidate optimal points obtained by the algorithm in both stochastic and deterministic settings to the known global optima for each benchmark function.

\begin{table}[h!]
\centering
\caption{Each row shows the benchmark function. The columns indicate the best solution under stochastic sampling, the best solution under purely deterministic sampling, and the known global optimum. The values are of the form $([x,y],f)$.}
\label{tab:opt}
\begin{tabular}{l l l l}
\toprule
\textbf{Function} 
& \textbf{Stochastic Solution} 
& \textbf{Deterministic Solution}
& \textbf{True Global Optimum} \\
\midrule
Rosenbrock  
& $[(0.71,\, 0.46],\, 0.21)$  
& $([0.98,\, 0.95],\, 0.00)$ 
& $([1.00,\, 1.00],\, 0.00)$ \\
Himmelblau 
& $[(2.88,\, 1.97],\, 0.61)$  
& $([3.59,\, -1.84],\, 0.01)$ 
& $([3.00,\, 2.00], [3.58,\, -1.85]\, 0.00)$,\\
Booth       
& $(0.71,\, 3.54],\, 0.61)$  
& $([0.99,\, 2.99],\, 0.00)$ 
& $([1.00,\, 3.00],\, 0.00)$ \\
Ackley      
& $(0.05,\, -0.05],\, 0.33)$ 
& $([-0.03,\, -0.03],\, 0.03)$ 
& $([0.00,\, 0.00],\, 0.00)$ \\
Rastrigin   
& $(-0.05,\, -0.05],\, 1.05)$ 
& $([-0.03,\, -0.03],\, 0.13)$ 
& $([0.00,\, 0.00],\, 0.00)$ \\
\bottomrule
\end{tabular}
\end{table}

As we can see from Table~\ref{tab:opt}, the deterministic optimal points lie very close to the actual global optimum, indicating that Algorithm~2 can indeed locate the global optimum when the surrogate precisely approximates the true function’s shape. Given the current tuning parameters, the deterministic optimal point is the best solution that Algorithm~2 can achieve. For Algorithm~2 to match these deterministic values, the surrogate manifold must closely approximate the actual function.

Additionally, the stochastic optimal points converge near the global optima. To further refine these points, one can iteratively resample the function around the stochastic optimum by reapplying Algorithm~1 with circles (geodesic hyper-sphere in high dimensions) closer to that region, then reapplying Algorithm~2. This process can be repeated until the desired level of accuracy is reached. Error mtrics can be used to make this decision. This hybrid approach is specially useful for high dimensional and/ or oscillatory functions. \\

To illustrate this, we apply one iteration of this hybrid approach, i.e. after getting optimal in the initial iteration, we recreate the manifold in the region close to the initial optima with the revised domain $[-1,1] \times [-1,1]$ with a Fourier series of order 5. Doing so improves the error metrics to to Mean Absolute Error (MAE): 0.1034, Mean Squared Error (MSE): 0.0688, R-squared: 0.9993. Additionally, the revised optimum is obtained at $([0.01, 0.010], 0.04)$ which is very close to true global optimum $([0.00, 0.00], 0)$ .\\

\begin{figure}[H]
   \hspace*{-2cm} 
    \includegraphics[width=1.2\linewidth]{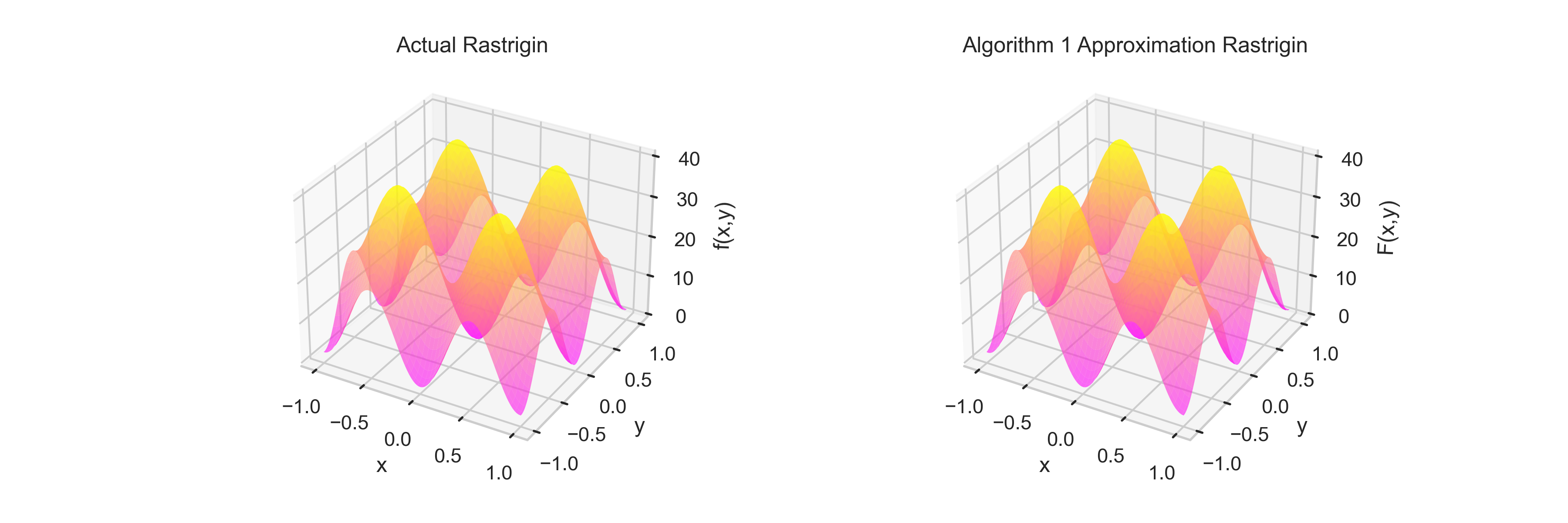}
    \caption{Rastrigin approximation using stochastic sampling in revised domain.}
     \label{Rastrigin}
\end{figure}

We can say that Algorithm~1 provides a surrogate manifold reasonably close to the standard global shape of the function. The algorithm, 2 successfully converges to candidate points that are close to the known global optima even for highly non-convex functions such as Rastrigin.

\section{Conclusion}
In this paper, we proposed a novel method for optimizing the stochastic processes. Using a defined sampling procedure that allows for boundary and center point sampling, a surrogate manifold is approximated using Fourier series approximation. The manifold is iteratively evolved using samples on the circles formed on the evolved matrix using the Riemannian metric tensor. We continue the procedure until the circle crosses the boundary of the manifold. We also establish error bounds on the surrogate manifold precision. 

For optimization, we regularize and de-regularize the manifold using a combination of Ricci flow and inverse Ricci flow such that the regions of large curvature tend to blow up to form singularities. These points form a set of candidates for optima. These points are checked and compared for the best solution. The solution obtained is the optimum.

We apply this algorithm to five benchmark functions: Rosenbrock, Himmeblau, Booth, Ackley and Rastrigin. We observed that Algorithm~1 can closely approximate the actual function shape for relatively smooth or flatter surfaces. For oscillatory functions, albeit approximating the local variations require close circle sampling, the global structure can be approximated qualitatively. 

Algorithm~2 finds the global optimum with high degree of accuracy given the original manifold is approximated perfectly. With stochastic sampling, the algorithm converged close to the true optima. However, to get more accurate solution, one can iteratively approximate the manifold near stochastic optima and find better solution. This procedure can be repeated until the solution stablizes.

\section{Data Availability Statement}
The manuscript has no associated data.
\section*{References}
\medskip

{
\small

[1] Powell, Warren B. "A unified framework for stochastic optimization." European Journal of Operational Research 275, no. 3 (2019): 795-821.

[2] Dantzig, George B. "Linear programming under uncertainty." Management science 1, no. 3-4 (1955): 197-206.

[3] Beale, Evelyn ML. "On minimizing a convex function subject to linear inequalities." Journal of the Royal Statistical Society Series B: Statistical Methodology 17, no. 2 (1955): 173-184.

[4] Shapiro, Alexander, Darinka Dentcheva, and Andrzej Ruszczynski. Lectures on stochastic programming: modeling and theory. Society for Industrial and Applied Mathematics, 2021.

[5] Li, C., \& Grossmann, I. E. (2021). A review of stochastic programming methods for optimization of process systems under uncertainty. Frontiers in Chemical Engineering, 2, 622241.

[6] Green, John W. "Approximately convex functions." (1952): 499-504.

[7] Lovász, László, and Santosh Vempala. "Fast algorithms for logconcave functions: Sampling, rounding, integration and optimization." In 2006 47th Annual IEEE Symposium on Foundations of Computer Science (FOCS'06), pp. 57-68. IEEE, 2006. 

[8] Kannan, Ravi, László Lovász, and Miklós Simonovits. "Random walks and an o*(n5) volume algorithm for convex bodies." Random Structures \& Algorithms 11, no. 1 (1997): 1-50.

[9] Lovász, László, and Miklós Simonovits. "Random walks in a convex body and an improved volume algorithm." Random structures \& algorithms 4, no. 4 (1993): 359-412.

[10] Lovász, László, and Santosh Vempala. "Hit-and-run from a corner." In Proceedings of the thirty-sixth annual ACM symposium on Theory of computing, pp. 310-314. 2004.

[11] Kalai, Adam Tauman, and Santosh Vempala. "Simulated annealing for convex optimization." Mathematics of Operations Research 31, no. 2 (2006): 253-266.

[12] Liang, Tengyuan, Hariharan Narayanan, and Alexander Rakhlin. "On zeroth-order stochastic convex optimization via random walks." arXiv preprint arXiv:1402.2667 (2014).

[13] Shamir, Ohad. "On the complexity of bandit and derivative-free stochastic convex optimization." In Conference on Learning Theory, pp. 3-24. PMLR, 2013.

[14] Homem-de-Mello, Tito, and Güzin Bayraksan. "Monte Carlo sampling-based methods for stochastic optimization." Surveys in Operations Research and Management Science 19, no. 1 (2014): 56-85.

[15] Aien, Morteza, Masoud Rashidinejad, and Mahmud Fotuhi-Firuzabad. "On possibilistic and probabilistic uncertainty assessment of power flow problem: A review and a new approach." Renewable and Sustainable energy reviews 37 (2014): 883-895.

[16] Chow, Bennett, Peng Lu, and Lei Ni. Hamilton's Ricci flow. Vol. 77. American Mathematical Soc., 2006.
[17] Hamilton, Richard S. "The formation of singularities in the Ricci flow." Surveys in Diff. Geom. 2 (1995): 7-136.

[18] Gu, Xianfeng, Yalin Wang, Tony F. Chan, Paul M. Thompson, and Shing-Tung Yau. "Genus zero surface conformal mapping and its application to brain surface mapping." IEEE transactions on medical imaging 23, no. 8 (2004): 949-958.

[19] Adams, Robert A., and John JF Fournier. Sobolev spaces. Vol. 140. Elsevier, 2003.

[20] Shapiro, Alexander, Darinka Dentcheva, and Andrzej Ruszczynski. Lectures on stochastic programming: modeling and theory. Society for Industrial and Applied Mathematics, 2021.

[21] Rosenbrock, HoHo. "An automatic method for finding the greatest or least value of a function." The computer journal 3, no. 3 (1960): 175-184.

[22] Himmelblau, David M. Applied nonlinear programming. McGraw-Hill, 2018.

[23] Surjanovic, Sonja, and Derek Bingham. Virtual library of simulation experiments: test functions and datasets. 2013.

[24] Ackley, David. A connectionist machine for genetic hillclimbing. Vol. 28. Springer science \& business media, 2012.

[25] Rastrigin, Leonard Andreevich. "About convergence of random search method in extremal control of multi‐parameter systems." Avtomat I Telemekh 24 (1963): 1467.

\end{document}